\documentclass[12pt]{article}
\usepackage{amsmath}
\usepackage{amssymb}
\usepackage{vatola}

\def\q{\quad}
\def\qq{\qquad}
\def\qtq#1{\q\t{#1}\q}
\def\mod{\pmod}
\def\t{\text}
\def\f{\frac}
\def\e{\equiv}
\def\b{\binom}
\def\dn#1{d_n^{(#1)}}
\def\sls#1#2{(\f{#1}{#2})}

\def\Ls#1#2{\Big(\f{#1}{#2}\Big)}
\let \pro=\proclaim
\let \endpro=\endproclaim

\begin{document}
\leftline{J. Math. Anal. Appl. 456(2017), no.2, 912-926}
\par\q\par\q
 \centerline {\bf A kind of orthogonal polynomials and related
 identities}

 \par\q\newline \centerline{Zhi-Hong Sun}
\par\q\newline \centerline{School of Mathematical Sciences}
\centerline{Huaiyin Normal University}
 \centerline{Huaian, Jiangsu 223300, P.R. China}
  \centerline{Email: zhsun@hytc.edu.cn}
\centerline{URL: http://www.hytc.edu.cn/xsjl/szh} \par\q
\newline
{\bf Abstract} In this paper we introduce the
 polynomials $\{d_n^{(r)}(x)\}$
  and $\{D_n^{(r)}(x)\}$ given by
$d_n^{(r)}(x)=\sum_{k=0}^n\b{x+r+k}k\b{x-r}{n-k}
 \ (n\ge 0)$,
$D_0^{(r)}(x)=1,\ D_1^{(r)}(x)=x$ and $D_{n+1}^{(r)}(x)
=xD_n^{(r)}(x)-n(n+2r)D_{n-1}^{(r)}(x)\ (n\ge 1).$
 We show that
$\{D_n^{(r)}(x)\}$
  are orthogonal polynomials for $r>-\f 12$, and
   establish many
identities for $\{d_n^{(r)}(x)\}$
  and $\{D_n^{(r)}(x)\}$, especially
  obtain a formula
  for $d_n^{(r)}(x)^2$ and the
  linearization formulas for
   $d_m^{(r)}(x)d_n^{(r)}(x)$ and
   $D_m^{(r)}(x)D_n^{(r)}(x)$. As an application we extend recent
work of Sun and Guo.

\par\q
\par\q

 \section*{1. Introduction}
 \par\q Let $\Bbb Z$, $\Bbb N_0$ and $\Bbb N$
  be the sets of integers, nonnegative integers and
  positive integers, respectively. By [5, (3.17)],
 for $n\in\Bbb N_0$,
$$\sum_{k=0}^n\b nk\b xk t^k
=\sum_{k=0}^n\b nk\b{x+k}n(t-1)^{n-k}.\tag 1.1$$ Define
$$d_n(x)=\sum_{k=0}^n\b nk\b xk 2^k\q(n=0,1,2,
\ldots).\tag 1.2$$ For $m,n\in\Bbb N$, $d_n(m)$ is the number of
lattice paths from $(0, 0)$ to $(m,n)$, with jumps $(0, 1),\ (1, 1)$
or $(1, 0)$. $\{d_n(m)\}$ are called Delannoy
 numbers. See [2].
  In [8] Z.W. Sun deduced some supercongruences involving
$d_n(x)$. Actually, he obtained congruences for
$$\sum_{k=0}^{p-1}d_k(x)^2,\ \sum_{k=0}^{p-1}(-1)^kd_k(x)^2
,\ \sum_{k=0}^{p-1}(2k+1)d_k(x)^2\qtq{and}
\sum_{k=0}^{p-1}(-1)^k(2k+1)d_k(x)^2\tag 1.3$$ modulo $p^2$, where
$p$ is an odd prime and $x$ is a rational $p$-adic integer.
  Z.W. Sun also conjectured that for any $n\in\Bbb N$ and
$x\in\Bbb Z$,
$$\align &x(x+1)\sum_{k=0}^{n-1}(2k+1)d_k(x)^2\e 0\mod{2n^2},
\tag 1.4
\\&\sum_{k=0}^{n-1}\varepsilon^k(2k+1)
d_k(x)^{2m}\e 0\mod n\q\t{for given $\varepsilon\in \{1,-1\}$ and
$m\in\Bbb N$}.\tag 1.5\endalign$$ Recently, Guo[6] proved the above
two congruences by using the identity
$$d_n(x)^2=\sum_{k=0}^n\b{n+k}{2k}\b xk
\b{x+k}k4^k.\tag 1.6$$ Guo proved (1.6) by using Maple and
Zeilberger's algorithm, and Zudilin stated that (1.6) can be deduced
from two transformation formulas for hypergeometric series. See [6]
and [7, (1.7.1.3) and (2.5.32)].
\par In this paper we establish closed formulas
 for sums in (1.3),
 which imply Sun's related
congruences. Set
$$d_n^{(r)}(x)=\sum_{k=0}^n\b{x+r+k}k\b{x-r}{n-k}\
(n=0,1,2,\ldots).\tag 1.7$$ Then $d_n(x)=d_n^{(0)}(x)$ by (1.1).
Thus, $d_n^{(r)}(x)$ is a generalization of $d_n(x)$. The main
purpose of this paper is to investigate the properties of
$d_n^{(r)}(x)$. We establish many identities for $d_n^{(r)}(x)$. In
particular, we obtain a formula for $d_n^{(r)}(x)^2$, which is a
generalization of (1.6). See Theorem 2.6.
\par Some classical orthogonal polynomials have formulas
for the linearization of their products. As examples, for Hermite
polynomials $\{H_n(x)\}$ ($H_{-1}(x)=0,\ H_0(x)=1,\
H_{n+1}(x)=2xH_n(x)-2nH_{n-1}(x)\ (n\ge 0)$) and
 Legendre
polynomials $\{P_n(x)\}$ ($P_0(x)=1,\ P_1(x)=x, \
(n+1)P_{n+1}(x)=(2n+1)xP_n(x)-nP_{n-1}(x)\ (n\ge 1)$) we have the
linearization of their products.
 See [1, Theorem 6.8.1 and
Corollary 6.8.3] and [3, p.195]. In Section 2 we establish the
following linearization formula:
$$d_m^{(r)}(x)d_n^{(r)}(x)=\sum_{k=0}^{\min\{m,n\}}\b{m+n-2k}{m-k}\b{2r+m+n-k}k
(-1)^kd_{m+n-2k}^{(r)}(x).\tag 1.8$$
 \par In Section 3 we introduce the polynomials
 $\{D_n^{(r)}(x)\}$ given by
$$D_0^{(r)}(x)=1,\  D_1^{(r)}(x)=x\qtq{and}  D_{n+1}^{(r)}(x)
=xD_n^{(r)}(x)-n(n+2r)D_{n-1}^{(r)}(x)\ (n\ge 1). \tag 1.9$$  By [4,
pp.175-176]
 or [1, pp.244-245], $\{D_n^{(r)}(x)\}$ are orthogonal
polynomials for $r>-\f 12$, although we have not found their weight
functions. We state that
$D_n^{(r)}(x)=(-i)^nn!d_n^{(r)}\sls{ix-1}2$, and obtain some
properties of $\{D_n^{(r)}(x)\}$. In particular, we show that
$$D_n^{(r)}(x)^2-D_{n+1}^{(r)}(x)D_{n-1}^{(r)}(x)>0\q
\t{for $r>-\f 12$ and real $x$}.\tag 1.10$$ Note that
$P_n(x)^2-P_{n-1}(x)P_{n+1}(x)\ge 0$ for $|x|\le 1$ and
$H_n(x)^2-H_{n-1}(x)H_{n+1}(x)\ge 0$. See [1, p.342] and
[3, p.195].
\par Throughout this paper, $[a]$ is the
greatest integer
 not exceeding $a$, and $f'(x)$ is the derivative of $f(x)$.

\section*{2. The properties of $d_n^{(r)}(x)$}

\par
By (1.1) and (1.2), for $n\in\Bbb N_0$,
$$d_n(x)=\sum_{k=0}^n\b nk\b xk2^k=\sum_{k=0}^n\b nk\b{x+k}n=\sum_{k=0}^n
\b {x+k}k\b x{n-k}.\tag 2.1$$
 Now we introduce the following generalization of $\{d_n(x)\}$.
 \pro{Definition 2.1} Let $\{d_n^{(r)}(x)\}$ be the polynomials
given by
$$d_n^{(r)}(x)=\sum_{k=0}^n\b{x+r+k}k\b{x-r}{n-k}\
(n=0,1,2,\ldots).$$ For convenience we also define
$d_{-1}^{(r)}(x)=0.$\endpro
\par By (2.1), $d_n(x)=d_n^{(0)}(x)$. Since $\b{-a}k=(-1)^k\b{a+k-1}k$ we see that
$$d_n^{(r)}(x)=\sum_{k=0}^n\b{-1-x-r}k(-1)^k\b{x-r}{n-k}
=\sum_{k=0}^n\b{-1-x-r}{n-k}(-1)^{n-k}\b{x-r}k.
\tag 2.2$$ Hence
$$d_n^{(r)}(-1-x)=(-1)^nd_n^{(r)}(x).\tag 2.3$$

  The first few $\{\dn r(x)\}$ are shown below:
$$\align &d_0^{(r)}(x)=1,\ d_1^{(r)}(x)=2x+1,\
d_2^{(r)}(x)=2x^2+2x+r+1,
\\&d_3^{(r)}(x)=\f 43x^3+2x^2+\big
(2r+\f 83\big)x+r+1.
\endalign$$

 \pro{Theorem 2.1} For $|t|<1$ we have
$$\sum_{n=0}^{\infty}d_n^{(r)}(x)t^n
=\f{(1+t)^{x-r}}{(1-t)^{x+r+1}}.\tag 2.4$$
\endpro
Proof. Newton's binomial theorem states that
$(1+t)^{\alpha}=\sum_{n=0}^{\infty}\b{\alpha}nt^n.$ Thus,
$$\align (1+t)^{x-r}(1-t)^{-x-r-1}&=\Big(\sum_{m=0}^{\infty}\b{x-r}mt^m\Big)\Big(\sum_{k=0}^{\infty}
\b{-x-r-1}k(-1)^kt^k\Big)
\\&=\sum_{n=0}^{\infty}\Big(\sum_{k=0}^n\b{-x-r-1}k(-1)^k
\b{x-r}{n-k}\Big)t^n =\sum_{n=0}^{\infty}d_n^{(r)}(x)t^n.\endalign$$
This proves the theorem.$\q\square$
 \pro{Corollary 2.1} For $n\in\Bbb N$ we have
$$d_n^{(r)}\Big(-\f 12\Big)=\cases 0&\t{if $2\nmid n$,}
\\\b{-1/2-r}{n/2}(-1)^{n/2}&\t{if $2\mid n$.}\endcases$$
\endpro
Proof. By Theorem 2.1 and Newton's binomial theorem, for $|t|<1$ we
have
$$\sum_{n=0}^{\infty}d_n^{(r)}(-1/2)t^n
=(1-t^2)^{-1/2-r}=\sum_{k=0}^{\infty}\b{-1/2-r}
k(-1)^kt^{2k}.$$ Now
comparing the coefficients
 of $t^n$ on both sides yields the result.
$\q\square$

\pro{Theorem 2.2} For $n\in\Bbb N$ we have
$$(n+1)d_{n+1}^{(r)}(x)=(1+2x)d_n^{(r)}(x)+(n+2r)
d_{n-1}^{(r)}(x).\tag 2.5$$
\endpro
Proof. By Theorem 2.1, for $|t|<1$,
$$\align &\sum_{n=0}^{\infty}(n+1)d_{n+1}^{(r)}(x)t^n
-\sum_{n=0}^{\infty}nd_{n-1}^{(r)}(x)t^n
\\&=\Big(\sum_{n=0}^{\infty}d_{n+1}^{(r)}(x)t^{n+1}\Big)'
-t\Big(\sum_{n=1}^{\infty}d_{n-1}^{(r)}(x)t^n\Big)'
\\&=\big((1+t)^{x-r}(1-t)^{-x-r-1}\big)'-t\big(t(1+t)^{x-r}(1-t)^{-x-r-1}\big)'
\\&=\big((1+t)^{x-r}(1-t)^{-x-r-1}\big)'
-t((1+t)^{x-r}(1-t)^{-x-r-1}+t((1+t)^{x-r}(1-t)^{-x-r-1})')
\\&=(1-t^2)\big((x-r)(1+t)^{x-r-1}(1-t)^{-x-r-1}+(1+t)^{x-r}(x+r+1)(1-t)^{-x-r-2}\big)
\\&\qq-t(1+t)^{x-r}(1-t)^{-x-r-1}
\\&=(1+2x+2rt)(1+t)^{x-r}(1-t)^{-x-r-1}
\\&=(1+2x)\sum_{n=0}^{\infty}d_n^{(r)}(x)t^n+2r\sum_{n=1}^{\infty}d_{n-1}^{(r)}(x)t^n.
\endalign$$
Now comparing the coefficients of $t^n$ on both sides gives the
result. $\q\square$

\pro{Theorem 2.3} Let $n\in\Bbb N_0$. Then
$$d_n^{(r)}(x)=\sum_{k=0}^{[n/2]}
\b{r-1+k}kd_{n-2k}(x)=\sum_{k=0}^n \b{2r-1+k}k d_{n-k}(x-r)\tag 2.6
$$
and
$$d_n(x)=\sum_{k=0}^{[n/2]}\b rk(-1)^kd_{n-2k}^{(r)}(x)
=\sum_{k=0}^n\b {2r}k(-1)^kd_{n-k}^{(r)}(x+r). \tag 2.7$$
\endpro
Proof. By (2.4),
$$\sum_{n=0}^{\infty}\dn r(x)t^n=(1-t^2)^{-r}\cdot\f
1{1-t}\Ls{1+t}{1-t}^x=(1-t)^{-2r}\cdot\f
1{1-t}\Ls{1+t}{1-t}^{x-r}.$$ Hence
$$\sum_{n=0}^{\infty}\dn r(x)t^n=(1-t^2)^{-r}
\sum_{n=0}^{\infty}
d_n(x)t^n=(1-t)^{-2r}\sum_{n=0}^{\infty}
 d_n(x-r)t^n,$$ which yields the first 2 results by
 applying
Newton's binomial theorem and comparing the coefficients of $t^n$ on
both sides.
 Also,
$$\sum_{n=0}^{\infty} d_n(x)t^n=(1-t^2)^r\sum_{n=0}^{\infty}\dn r(x)t^n
=(1-t)^{2r}\sum_{n=0}^{\infty}\dn r(x+r)t^n$$ yields the next 2
results. $\q\square$

\pro{Corollary 2.2} Let $n\in\Bbb N_0$. Then $d_n^{(r)}(0)=\b{r+[\f
n2]}{[\f n2]}.$
\endpro
Proof. Set $\b ak=0$ for $k<0$. Since $d_n(0)=\sum_{k=0}^n\b nk\b
0k2^k=1$, applying Theorem 2.3 we get
$$\align d_n^{(r)}(0)&=\sum_{k=0}^{[n/2]}\b{r-1+k}k
=\sum_{k=0}^{[n/2]}\b{-r}k(-1)^k
\\&=\sum_{k=0}^{[n/2]}\Big((-1)^k\b{-r-1}k-(-1)^{k-1}\b{-r-1}{k-1}
\Big) =(-1)^{[\f n2]}\b{-r-1}{[\f n2]}=\b{r+[\f n2]}{[\f
n2]}.\q\square\endalign$$

\pro{Theorem 2.4} For $n\in\Bbb N$ we have
$$\align &d_n^{(r)}(x)=d_n^{(r+1)}(x)-d_{n-2}^{(r+1)}(x),\tag i
\\&d_n^{(r+1)}(x)
=\sum_{k=0}^{[n/2]}d_{n-2k}^{(r)}(x),\tag ii
\\&(n+1)^2d_{n+1}^{(r)}(x)^2-(n+2r+1)^2d_n^{(r)}(x)^2
=4(x-r)(x+1+r)(d_n^{(r+1)}(x)^2-d_{n-1}^{(r+1)}(x)^2),\tag iii
\\&(2r+1)\sum_{k=0}^{n-1}(2k+2r+1)d_k^{(r)}(x)^2
=n^2d_n^{(r)}(x)^2-4(x-r)(x+1+r)d_{n-1}^{(r+1)}(x)^2.\tag iv
\endalign$$
\endpro
Proof. By Theorem 2.1, for $|t|<1$,
$$\sum_{n=0}^{\infty}d_n^{(r+1)}(x)t^n=\f
1{1-t^2}\sum_{m=0}^{\infty}d_m^{(r)}(x)t^m=
\Big(\sum_{k=0}^{\infty}t^{2k}\Big)
\Big(\sum_{m=0}^{\infty}d_m^{(r)}(x)t^m\Big).$$ Now comparing the
coefficients of $t^n$ on both sides yields (i) and (ii).
\par By (i) and (2.5),
$$\align &d_{n+1}^{(r)}(x)=d_{n+1}^{(r+1)}(x)-d_{n-1}^{(r+1)}(x)
\\&=\f{(2x+1)d_n^{(r+1)}(x)+(n+2+2r)d_{n-1}^{(r+1)}(x)}{n+1}
-d_{n-1}^{(r+1)}(x)
=\f{(2x+1)d_n^{(r+1)}(x)+(2r+1)d_{n-1}^{(r+1)}(x)}{n+1}\endalign$$
and $$\align &d_n^{(r)}(x)=d_n^{(r+1)}(x) -d_{n-2}^{(r+1)}(x)
\\&=d_n^{(r+1)}(x)-\f{nd_n^{(r+1)}(x)-(2x+1)d_{n-1}^{(r+1)}(x)}{n+1+2r}
=\f{(2r+1)d_n^{(r+1)}(x)+(2x+1)d_{n-1}^{(r+1)}(x)}{n+1+2r}.\endalign$$
Thus,
$$\align &(n+1)^2d_{n+1}^{(r)}(x)^2-(n+1+2r)^2d_n^{(r)}(x)^2
\\&=((2x+1)d_n^{(r+1)}(x)+(2r+1)d_{n-1}^{(r+1)}(x))^2
-((2r+1)d_n^{(r+1)}(x)+(2x+1)d_{n-1}^{(r+1)}(x))^2
\\&=4(x-r)(x+1+r)(d_n^{(r+1)}(x)^2-d_{n-1}^{(r+1)}(x)^2).\endalign$$
This proves (iii). By (iii),
$$\align &\sum_{k=0}^{n-1}(2r+1)(2k+2r+1)d_k^{(r)}(x)^2
\\&=\sum_{k=0}^{n-1}\big((k+1)^2d_{k+1}^{(r)}(x)^2-k^2d_k^{(r)}(x)^2\big)
-4(x-r)(x+1+r)\sum_{k=0}^{n-1}\big(d_k^{(r+1)}(x)^2
-d_{k-1}^{(r+1)}(x)^2\big)
\\&=n^2d_n^{(r)}(x)^2-4(x-r)(x+1+r)d_{n-1}^{(r+1)}(x)^2.\endalign$$
This proves (iv). $\q\square$

\pro{Theorem 2.5} Let $n\in\Bbb N$, $r\in\Bbb N_0$ and $x\in\Bbb Z$.
Then
$$(2r+1)\prod_{k=-r}^r(x+k)(x+1-k)
\sum_{k=0}^{n-1}(2k+2r+1)d_k^{(r)}(x)^2
\e 0\mod{2n^2(n+1)^2\cdots
(n+2r)^2}.$$
\endpro
Proof. It is easily seen that for $k,n,r\in\Bbb N_0$
 with $k\le n$,
$$\b{x+r}{2r}\b{x+r+k}k\b{x-r}{n-k}=\b{n+2r}{2r}\b
nk\b{x+r+k}{n+2r}.$$ Thus,
$$\b{x+r}{2r}d_n^{(r)}(x)=\b{n+2r}{2r}
\sum_{k=0}^n\b nk\b{x+r+k}{n+2r} \qtq{for}r\in\Bbb N_0. \tag 2.8$$
By Theorem 2.4(iv) and (2.8),
$$\align &(2r+1)\prod_{k=-r}^r(x+k)(x+1-k)
\sum_{k=0}^{n-1}(2k+2r+1)d_k^{(r)}(x)^2
\\&=\prod_{k=-r}^r(x+k)(x+1-k)\times \big(n^2d_n^{(r)}(x)^2-4(x-r)(x+1+r)
d_{n-1}^{(r+1)}(x)^2\big)
\\&=(x-r)(x+r+1)(n+2r)^2(n+2r-1)^2\cdots(n+1)^2n^2\Big(\sum_{k=0}^n\b nk\b{x+r+k}{n+2r}\Big)^2
\\&\q-4(n+2r+1)^2(n+2r)^2\cdots n^2
\Big(\sum_{k=0}^{n-1}\b{n-1}k \b{x+r+1+k}{n+2r+1}\Big)^2.
\endalign$$
To finish the proof, we note that
 $(x+r+1)(x-r)\e 0\mod 2$. $\q\square$
\par We remark that Theorem 2.5
 is a generalization of (1.4),
and the next theorem is a generalization of
 (1.6).

\pro{Theorem 2.6} Suppose $n\in\Bbb N_0$ and $r\not\in\{-\f 12,-\f
22,-\f 32,\ldots\}$. Then
$$d_{n}^{(r)}(x)^2
=\b{n+2r}{n}\sum_{m=0}^n\f{\b{x-r}m \b{x+r+m}m\b{n+2r+m}{n-m}}
{\b{m+2r}{m}}4^m.\tag 2.9$$
\endpro
Proof. Set
$$s(n)=\f{d_n^{(r)}(x)^2}{\b{n+2r}n}\qtq{and}
S(n)=\sum_{m=0}^n\f{\b{x-r}m \b{x+r+m}m\b{n+2r+m}{n-m}}
{\b{m+2r}{m}}4^m.$$ Using sumrecursion in Maple we find that for
$n\in\Bbb N$,
$$(n+2)(n+2+2r)S(n+2)-((2x+1)^2+(n+1)(n+1+2r))
(S(n+1)+S(n))+n(n+2r)S(n-1)=0.$$ By Theorem 2.2,
$$d_{n+2}^{(r)}(x)=\f{(1+2x)d_{n+1}^{(r)}(x)
+(n+1+2r)d_n^{(r)}(x)}{n+2},\
d_{n-1}^{(r)}(x)=\f{(n+1)d_{n+1}^{(r)}(x)
-(1+2x)d_n^{(r)}(x)}{n+2r}.$$ Thus,
$$\align&(n+2)(n+2+2r)s(n+2)+n(n+2r)s(n-1)
\\&=\f{(n+2)(n+2+2r)}{\b{n+2+2r}{2r}}d_{n+2}^{(r)}(x)^2
+\f{n(n+2r)}{\b{n-1+2r}{2r}}d_{n-1}^{(r)}(x)^2
\\&=\f {((1+2x)d_{n+1}^{(r)}(x)+(n+1+2r)d_n^{(r)}(x))^2}{\b{n+1+2r}{2r}}
+\f{((n+1)d_{n+1}^{(r)}(x)-(1+2x)d_n^{(r)}(x))^2}{\b{n+2r}{2r}}
\\&=d_{n+1}^{(r)}(x)^2\Big\{\f{(1+2x)^2}{\b{n+1+2r}{2r}}+\f{(n+1)^2}
{\b{n+2r}{2r}}\Big\} +
d_n^{(r)}(x)^2\Big\{\f{(1+2x)^2}{\b{n+2r}{2r}}+\f{(n+1+2r)^2}
{\b{n+1+2r}{2r}}\Big\}
\\&=\f{d_{n+1}^{(r)}(x)^2}{\b{n+1+2r}{2r}}\big((1+2x)^2+(n+1)(n+1+2r)\big)
+\f{d_n^{(r)}(x)^2}{\b{n+2r}{2r}} \big((1+2x)^2+(n+1)(n+1+2r)\big) .
\\&=((1+2x)^2+(n+1)(n+1+2r))(s(n)+s(n+1)).
\endalign$$
This shows that $s(n)$ and $S(n)$ satisfy the same recurrence
relation. Also,
$$\align &s(0)=1=S(0),\ s(1)=\f{(1+2x)^2}{2r+1}=S(1),
\  s(2)=\f{(2x^2+2x+r+1)^2}{(r+1)(2r+1)}=S(2).
\endalign$$
Thus, $s(n)=S(n)$ for $n\in\Bbb N_0$. $\q\square$
\par
\q \par Now we present the linearization of
$d_m^{(r)}(x)d_n^{(r)}(x)$.
 \pro{Theorem 2.7} Let $m,n\in\Bbb N_0$. Then
$$d_m^{(r)}(x)d_n^{(r)}(x)=\sum_{k=0}^{\min\{m,n\}}\b{m+n-2k}{m-k}\b{2r+m+n-k}k
(-1)^kd_{m+n-2k}^{(r)}(x).\tag 2.10$$
\endpro
\par Proof. Let $L(m,n)=d_m^{(r)}(x)d_n^{(r)}(x)$
 and $\b ak=0$ for $k<0$.
By Theorem 2.2, $(m+1+2r)d_m^{(r)}(x)+(1+2x)d_{m+1}^{(r)}(x) =(m+2)
d_{m+2}^{(r)}(x)$. Hence
$$(m+1+2r)L(m,n)+(1+2x)L(m+1,n)-(m+2)L(m+2,n)=0.$$
 Let
$$G(m,n,k,l) = (-1)^k {m+n-2k \choose m-k} {2r+m+n-k \choose k}
{x+r+l \choose l} {x-r \choose m+n-2k-l}.$$
 Using Maple it is easy
to check that
$$\align &(m+1+2r)G(m,n,k,l) + (2x+1)G(m+1,n,k,l)-(m+2)G(m+2,n,k,l) \\
\\&=F_1(m,n,k+1,l)-F_1(m,n,k,l)+F_2(m,n,k,l+1)-F_2(m,n,k,l),
\endalign$$
where
$$\align F_1(m,n,k,l)&=(-1)^k(2m+n+2r+4-2k)
\\&\q\times\b{m+n+2-2k}{m+2-k} \b{2r+m+1+n-k}{k-1}\b{x+r+l}l\b{x-r}{m+2+n-2k-l}
\endalign$$  and

$$\align &F_2(m,n,k,l) \\&=(-1)^kl
\b{m+1+n-2k}{m+1-k} \b{2r+m+n+1-k}k \b{x+r+l}l\b{x+1-r}{m+2+n-2k-l}
.\endalign$$ Thus,
$$\align&\sum_{k=0}^{m+2}\sum_{l=0}^{m+2+n}\Big((m+1+2r)G(m,n,k,l)
+(2x+1)G(m+1,n,k,l) -(m+2)G(m+2,n,k,l)\Big)
\\&=\sum_{l=0}^{m+2+n}\sum_{k=0}^{m+2}\big(F_1(m,n,k+1,l)-F_1(m,n,k,l)\big)
+\sum_{k=0}^{m+2}\sum_{l=0}^{m+2+n}\big(F_2(m,n,k,l+1)
-F_2(m,n,k,l)\big)
\\&=\sum_{l=0}^{m+2+n}(F_1(m,n,m+3,l)-F_1(m,n,0,l))+\sum_{k=0}^{m+2}
(F_2(m,n,k,m+n+3)-F_2(m,n,k,0))\\&=0.
\endalign$$
Set $$R(m,n) = \sum_{k=0}^m \sum_{l=0}^{m+n} G(m,n,k,l) =
\sum_{k=0}^m \sum_{l=0}^{m+n-2k} G(m,n,k,l).$$ Then
$(m+1+2r)R(m,n)+(2x+1)R(m+1,n)-(m+2)R(m+2,n)=0.$ From the above we
see that $L(m,n)$ and $R(m,n)$ satisfy the same recurrence relation.
It is clear that $L(0,n)=d_n^{(r)}(x)= \sum_{l=0}^n {x+r+l \choose
l}{x-r \choose n-l} = R(0,n).$ By Theorem 2.2,
$R(1,n)=(n+1)d_{n+1}^{(r)}(x)
-(n+2r)d_{n-1}^{(r)}(x)=(1+2x)d_n^{(r)}(x)=L(1,n).$ Hence, $L(m,n) =
R(m,n)$ for any nonnegative integers $m$ and $n$. This proves the
theorem. $\q\square$

\pro{Theorem 2.8} For $n\in\Bbb N$
 we have
$$\aligned &2(1+x+y)\sum_{k=0}^{n-1}\f
{(2r+k+1)\cdots (2r+n)}{(k+1)\cdots n}
 d_k^{(r)}(x)d_k^{(r)}(y)
\\&=(n+2r)(d_n^{(r)}(x)d_{n-1}^{(r)}(y)+d_{n-1}^{(r)}(x)
d_n^{(r)}(y)).\endaligned\tag 2.11$$
\endpro
Proof. We prove (2.11) by induction on $n$. Clearly (2.11)
is true
for $n=1$. By Theorem 2.2,
$$\align &(n+1)\big(d_{n+1}^{(r)}(x)d_n^{(r)}(y)+d_n^{(r)}(x)
d_{n+1}^{(r)}(y)\big)
\\&=d_n^{(r)}(y)\big((1+2x)d_n^{(r)}(x)+(n+2r)d_{n-1}^{(r)}(x)
\big)+d_n^{(r)}(x)\big((1+2y)d_n^{(r)}(y)+(n+2r)d_{n-1}^{(r)}(y)
\big)
\\&=2(1+x+y)d_n^{(r)}(x)d_n^{(r)}(y)+(n+2r)\big(d_n^{(r)}(x)d_{n-1}^{(r)}(y)+d_{n-1}^{(r)}(x)
d_n^{(r)}(y)\big).
\endalign$$ Thus,
if the result holds for $n$, then
$$\align &2(1+x+y)\sum_{k=0}^n\f{(2r+k+1)\cdots(2r+n+1)}
{(k+1)\cdots(n+1)}
d_k^{(r)}(x)d_k^{(r)}(y)
\\&=\f{n+2r+1}{n+1}2(1+x+y)
\Big(d_n^{(r)}(x)d_n^{(r)}(y)+\sum_{k=0}^{n-1}\f {(2r+k+1)\cdots
(2r+n)}{(k+1)\cdots n} d_k^{(r)}(x)d_k^{(r)}(y)\Big)
\\&=\f{n+2r+1}{n+1}\big(2(1+x+y)d_n^{(r)}(x)d_n^{(r)}(y)
+(n+2r)(d_n^{(r)}(x)d_{n-1}^{(r)}(y) +d_{n-1}^{(r)}(x)
d_n^{(r)}(y))\big)
\\&=(n+1+2r)\big(d_{n+1}^{(r)}(x)d_n^{(r)}(y)+d_n^{(r)}(x)
d_{n+1}^{(r)}(y)\big).
\endalign$$
Hence (2.11) holds for $n+1$. $\q\square$

\par{\bf Remark 2.1.}
Taking $r=0$ in Theorem 2.8 and noting that $d_n(x)=d_n^{(0)}(x)$
yields
$$2(1+x+y)\sum_{k=0}^{n-1}d_k(x)d_k(y)
=n(d_n(x)d_{n-1}(y)+d_{n-1}(x)d_n(y)).\tag 2.12$$

\section*{3. The orthogonal polynomials  $\{D_n^{(r)}(x)\}$}

\par By [4, pp.175-176], every orthogonal system of real valued
polynomials $\{p_n(x)\}$ satisfies
$$p_{-1}(x)=0,\ p_0(x)=1\qtq{and} xp_n(x)=A_np_{n+1}(x)+B_np_n(x)
+C_np_{n-1}(x)\ (n\ge 0),\tag 3.1$$ where $A_n,B_n,C_n$ are real and
$A_nC_{n+1}>0$. Conversely, if (3.1) holds for a sequence of
polynomials $\{p_n(x)\}$ and $A_n,B_n,C_n$ are real with
$A_nC_{n+1}>0$, then there exists a weight function $w(x)$ such that
$$\int_{-\infty}^{\infty}w(x)p_m(x)p_n(x)dx=\cases
0&\t{if $m\not=n$,}
\\\f 1{v_n}\int_{-\infty}^{\infty}w(x)dx&\t{if $m=n$,}
\endcases$$
where $v_0=1$ and $v_n=\f{A_0A_1\cdots A_{n-1}}{C_1\cdots C_n}$
 $(n\ge 1)$.

\par In this section we discuss a kind of orthogonal polynomials
related to $\{d_n^{(r)}(x)\}$.

 \pro{Definition 3.1} Let
$\{D_n^{(r)}(x)\}$ be the polynomials given by
$$D_{-1}^{(r)}(x)=0,\  D_0^{(r)}(x)=1\qtq{and}  D_{n+1}^{(r)}(x)
=xD_n^{(r)}(x)-n(n+2r)D_{n-1}^{(r)}(x)\ (n\ge 0). \tag 3.2$$\endpro
 \par The
first few $D_n^{(r)}(x)$ are shown below:
$$D_0^{(r)}(x)=1,\  D_1^{(r)}(x)=x,
\ D_2^{(r)}(x)=x^2-2r-1,\
  D_3^{(r)}(x)=x^3-(6r+5)x.$$
\par Suppose $r>-\f 12$. Set $A_n=1,\ B_n=0,\ C_n=n(n+2r)$, $v_0=1$ and $v_n=\f
1{n!(2r+1)(2r+2)\cdots(2r+n)}$ $(n\ge 1)$. Then $A_nC_{n+1}>0$ and
(3.1) holds for $p_n(x)=D_n^{(r)}(x)$. Hence $\{D_n^{(r)}(x)\}$ are
orthogonal polynomials.

\pro{Lemma 3.1} For $n\in\Bbb N_0$ we have
 $$d_n^{(r)}(x)=\f{i^nD_n^{(r)}(-i(1+2x))}{n!}\qtq{and so}
D_n^{(r)}(x)=(-i)^nn!d_n^{(r)}\Ls{ix-1}2.\tag 3.3$$
 \endpro
 Proof. Since $D_0^{(r)}(-i(1+2x))=1$,
 $iD_1^{(r)}(-i(1+2x))
 =1+2x$ and
 $$\align &(n+1)\f{i^{n+1}D_{n+1}^{(r)}(-i(1+2x))}
 {(n+1)!}
 \\&=\f{i^{n+1}D_{n+1}^{(r)}(-i(1+2x))}{n!}
 = \f{i^{n+1}}{n!}\big(-i(1+2x)D_n^{(r)}(-i(1+2x))-n(n+2r)D_{n-1}^{(r)}(-i(1+2x))\big)
\\&=(1+2x)\f{i^nD_n^{(r)}(-i(1+2x))}{n!}
+(n+2r)\f{i^{n-1}D_{n-1}^{(r)}(-i(1+2x))}
 {(n-1)!},\endalign$$ we must have
$d_n^{(r)}(x)=\f{i^nD_n^{(r)}(-i(1+2x))}{n!}$ by (2.5). Substituting
$x$ with $\f{ix-1}2$ yields the remaining part. $\q\square$
 \pro{Theorem 3.1} For $n\in\Bbb N$ we have
$$\sum_{k=0}^{n-1}(2k+2r+1)\prod_{s=k+1}^ns(s+2r)
 D_k^{(r)}(x)^2 =n(n+2r)\big(D_n^{(r)}(x)^2
 -D_{n-1}^{(r)}(x)D_{n+1}^{(r)}(x)\big).\tag 3.4$$
Thus, $D_n^{(r)}(x)^2-D_{n+1}^{(r)}(x)D_{n-1}^{(r)}(x)>0 $ for
$r>-\f 12$ and real $x$.
\endpro
Proof. Set
$\Delta_n^{(r)}(x)=D_n^{(r)}(x)^2-D_{n+1}^{(r)}(x)D_{n-1}^{(r)}(x)$.
We prove (3.4) by induction on $n$. Clearly (3.4) is true for $n=1$.
Suppose that (3.4) holds for $n$. Since
$$\align \Delta_{n+1}^{(r)}(x)
-n(n+2r)\Delta_n^{(r)}(x)
&=D_{n+1}^{(r)}(x)^2-D_n^{(r)}(x)(xD_{n+1}^{(r)}
(x)-(n+1)(n+2r+1)D_n^{(r)}(x))
\\&\q-n(n+2r)(D_n^{(r)}(x)^2-D_{n-1}^{(r)}(x)D_{n+1}^{(r)}(x))
\\&=D_{n+1}^{(r)}(x)(D_{n+1}^{(r)}(x)-
xD_n^{(r)}(x)+n(n+2r)D_{n-1}^{(r)}(x)) \\&\q
+((n+1)(n+1+2r)-n(n+2r))D_n^{(r)}(x)^2
\\&=(2n+2r+1)D_n^{(r)}(x)^2,
\endalign$$
we see that
$$\align &\sum_{k=0}^n(2k+2r+1)\prod_{s=k+1}^{n+1}s(s+2r)
\times D_k^{(r)}(x)^2
\\&=(n+1)(n+1+2r)\Big((2n+2r+1)D_n^{(r)}(x)^2+\sum_{k=0}^{n-1}(2k+2r+1)
\prod_{s=k+1}^ns(s+2r) D_k^{(r)}(x)^2\Big)
\\&=(n+1)(n+1+2r)\big((2n+2r+1)D_n^{(r)}(x)^2
+n(n+2r)\Delta_n^{(r)}(x)\big)
\\&=(n+1)(n+1+2r)\Delta_{n+1}^{(r)}(x).
\endalign$$
This shows that (3.4) holds for $n+1$. Hence (3.4) is proved by
induction. For $r>-\f 12$ we have
 $1+2r>0$. From (3.4) and the fact
 $D_0^{(r)}(x)=1$ we deduce that $\Delta_n^{(r)}(x)
 \ge (2r+1)\f{n!(2r+1)\cdots(2r+n)}
 {n(n+2r)}>0$. This concludes the proof.
$\q\square$

 \pro{Corollary 3.1} Let $n\in\Bbb N$. Then
$$\aligned &\sum_{k=0}^{n-1}(-1)^k(2k+2r+1)\f{(k+1+2r)\cdots(n+2r)}{(k+1)
\cdots n}d_k^{(r)}(x)^2 \\&=(-1)^n(n+2r)\big(nd_n^{(r)}(x)^2-(n+1)
d_{n-1}^{(r)}(x)d_{n+1}^{(r)}(x)\big).\endaligned\tag 3.5$$
\endpro
Proof. Replacing $x$ with $-i(1+2x)$ in Theorem 3.1 and then
applying Lemma 3.1 yields the result. $\q\square$

\pro{Theorem 3.2} Let $n\in\Bbb N$. Then
$$\sum_{k=0}^{n-1}\prod_{s=k+1}^ns(s+2r) D_k^{(r)}(x)^2
=n(n+2r)\Big(D_{n-1}^{(r)}(x)\f d{dx}D_n^{(r)}(x) -D_n^{(r)}(x)\f
d{dx}D_{n-1}^{(r)}(x)\Big)\tag 3.6$$ and
$$\aligned&\sum_{k=0}^{n-1}(-1)^k\prod_{s=k+1}^n\f{s+2r}s
\  d_k^{(r)}(x)^2 \\&=(-1)^{n-1}\f{n+2r}2\Big(d_{n-1}^{(r)}(x)\f
d{dx}d_n^{(r)}(x)-d_n^{(r)}(x)\f d{dx}
 d_{n-1}^{(r)}(x)\Big).\endaligned\tag
3.7$$
\endpro
Proof. We prove (3.6) by induction on $n$. Clearly (3.6)
 is true for
$n=1$. Suppose that (3.6) holds for $n$. Since
$D_{n+1}^{(r)}(x)=xD_n^{(r)}(x)-n(n+2r)D_{n-1}^{(r)}(x)$ we see that
$$\f d{dx}D_{n+1}^{(r)}(x)=D_n^{(r)}(x)+x\f
d{dx}D_n^{(r)}(x)-n(n+2r)\f d{dx}D_{n-1}^{(r)}(x)$$ and so
$$\align &D_n^{(r)}(x)\f d{dx}D_{n+1}^{(r)}(x)
-D_{n+1}^{(r)}(x) \f d{dx}D_n^{(r)}(x) -n(n+2r)
\Big(D_{n-1}^{(r)}(x)\f d{dx}D_n^{(r)}(x)-D_n^{(r)}(x)\f
d{dx}D_{n-1}^{(r)}(x)\Big)
\\&=D_n^{(r)}(x)^2+xD_n^{(r)}(x)
\f d{dx}D_n^{(r)}(x)-n(n+2r)D_n^{(r)}(x) \f d{dx}D_{n-1}^{(r)}(x)
\\&\q-(xD_n^{(r)}(x)-n(n+2r)D_{n-1}^{(r)}(x))\f d{dx}D_n^{(r)}(x)
\\&\q-n(n+2r)D_{n-1}^{(r)}(x)\f d{dx}D_n^{(r)}(x)
+n(n+2r)D_n^{(r)}(x)\f d{dx}D_{n-1}^{(r)}(x)
\\&=D_n^{(r)}(x)^2.\endalign$$
Hence
$$\align &\sum_{k=0}^n\prod_{s=k+1}^{n+1}s(s+2r)
\times D_k^{(r)}(x)^2
\\&=(n+1)(n+1+2r)
\Big(D_n^{(r)}(x)^2+\sum_{k=0}^{n-1}\prod_{s=k+1}^ns(s+2r)\cdot
D_k^{(r)}(x)^2\Big)
\\&=(n+1)(n+1+2r)
\Big(D_n^{(r)}(x)^2+ n(n+2r)\Big(D_{n-1}^{(r)}(x)\f
d{dx}D_n^{(r)}(x) -D_n^{(r)}(x)\f d{dx}D_{n-1}^{(r)}(x)\Big)\Big)
\\&=(n+1)(n+1+2r) \Big( D_n^{(r)}(x)\f d{dx}D_{n+1}^{(r)}(x)
-D_{n+1}^{(r)}(x) \f d{dx}D_n^{(r)}(x)\Big).
\endalign$$
This shows that (3.6) holds for $n+1$. Hence
 (3.6) is proved.

\par By Lemma 3.1, $d_n^{(r)}(x)=i^nD_n^{(r)}(-i(1+2x))/n!$. Thus,
$\f d{dx}d_n^{(r)}(x) =i^n\f d{dx}D_n^{(r)}(-i(1+2x))(-2i)/n!$.
 Now
applying (3.6) we obtain
$$\align &\sum_{k=0}^{n-1}(-1)^k
\prod_{s=k+1}^n\f{s+2r}s \times d_k^{(r)}(x)^2
\\&=\sum_{k=0}^{n-1}\prod_{s=k+1}^n\f{s+2r}s
\times\f{D_k^{(r)}(-i(1+2x))^2}{k!^2} =\f
1{n!^2}\sum_{k=0}^{n-1}\prod_{s=k+1}^n s(s+2r) \times
D_k^{(r)}(-i(1+2x))^2
\\&=\f{n(n+2r)}{n!^2}\Big(D_{n-1}^{(r)}(-i(1+2x))
\f d{dx}D_n^{(r)}(-i(1+2x)) -D_n^{(r)}(-i(1+2x))\f
d{dx}D_{n-1}^{(r)}(-i(1+2x))\Big)
 \\&=\f{n(n+2r)}{n!^2}\Big(\f{n!\f
d{dx}d_n^{(r)}(x)}{(-2i)i^n}\times
\f{(n-1)!d_{n-1}^{(r)}(x)}{i^{n-1}} - \f{n!d_n^{(r)}(x)}{i^n}\times
\f{(n-1)!\f d{dx}d_{n-1}^{(r)}(x)}{(-2i)i^{n-1}}\Big)
\\&=(-1)^{n-1}\f{n+2r}2\big(d_{n-1}^{(r)}(x)\f d{dx}d_n^{(r)}(x)-d_n^{(r)}(x)\f
d{dx}d_{n-1}^{(r)}(x)\big).
\endalign$$
This proves (3.7). $\q\square$

\par{\bf Remark 3.1.}  Taking $r=0$ in
(3.7) and (3.5) yields
$$\align &\sum_{k=0}^{n-1}(-1)^kd_k(x)^2=(-1)^{n-1}\f{n}2
(d_{n-1}(x)d_n'(x)-d_n(x)d_{n-1}'(x)), \tag 3.8
\\&\sum_{k=0}^{n-1}(-1)^k(2k+1)d_k(x)^2=(-1)^n(n^2d_n(x)^2-n(n+1)
d_{n-1}(x)d_{n+1}(x)).\tag 3.9\endalign$$

 \pro{Theorem 3.3} For $n\in\Bbb N_0$ we have
$$ D_n^{(r)}(x)^2=\sum_{m=0}^n\b{n+2r+m}{n-m}(-1)^{n-m}
\prod_{j=m+1}^nj(2r+j) \prod_{k=1}^m(x^2+(2r+2k-1)^2).\tag 3.10$$
\endpro
Proof. By Lemma 3.1 and Theorem 2.6,
$$\align D_n^{(r)}(x)^2=(-1)^nn!^2d_n^{(r)}\Ls{ix-1}2^2
=(-1)^nn!^2\b{n+2r}{n}\sum_{m=0}^n\f{\b{\f{ix-1}2-r}m
\b{\f{ix-1}2+r+m}m\b{n+2r+m}{n-m}} {\b{m+2r}{m}}4^m. \endalign$$
Since
$$\align &\b{\f{ix-1}2-r}m\b{\f{ix-1}2+r+m}m
\\&=\f{(\f{ix-1}2-r)(\f{ix-1}2-(r+1))\cdots(
\f{ix-1}2-(r+m-1))(\f{ix-1}2+r+m)\cdots(\f{ix-1}2+r+1)}{m!^2}
\\&=\f{((ix)^2-(2r+1)^2)\cdots((ix)^2-(2r+2m-1)^2)}{2^{2m}\cdot
m!^2}=\f{(x^2+(2r+1)^2)\cdots(x^2+(2r+2m-1)^2)}{(-4)^m\cdot m!^2}
,\endalign$$ from the above we deduce that
$$ D_n^{(r)}(x)^2=(-1)^nn!
\sum_{m=0}^n\b{n+2r+m}{n-m}\f{(-1)^m(2r+1)(2r+2)\cdots(2r+n)}
{m!(2r+1)(2r+2)\cdots(2r+m)} \prod_{k=1}^m(x^2+(2r+2k-1)^2).$$ This
yields the result. $\q\square$
 \pro{Theorem 3.4} The exponential
generating function of
 $\{D_n^{(r)}(x)\}$ is given by
$$\sum_{n=0}^{\infty}D_n^{(r)}(x)\f{t^n}{n!}=(1+t^2)^{-r-\f 12}
\;\t{\rm e}^{x\;arctan\; t}.\tag 3.11$$
\endpro
Proof.  Set $f(t)=\sum_{n=0}^{\infty}\  D_n^{(r)}(x)\f{t^n}{n!}$.
Then
$$f(t)=1+\sum_{n=0}^{\infty} D_{n+1}^{(r)}(x)\f{t^{n+1}}{(n+1)!}
=1+\sum_{n=0}^{\infty}x D_n^{(r)}(x)\f{t^{n+1}}{(n+1)!}
-\sum_{n=1}^{\infty}n(n+2r) D_{n-1}^{(r)}(x)\f{t^{n+1}}{(n+1)!}.$$
Hence
$$\align f'(t)&=\sum_{n=0}^{\infty}x D_n^{(r)}(x)\f{t^n}{n!}
-\sum_{n=1}^{\infty}(n+2r) D_{n-1}^{(r)}(x)\f{t^n}{(n-1)!}
\\&=xf(t)-2rtf(t)-t\Big(\sum_{n=1}^{\infty} D_{n-1}^{(r)}(x)
\f{t^n}{(n-1)!}\Big)'
\\&=(x-2rt)f(t)-t(tf(t))'=(x-2rt)f(t)-t(f(t)+tf'(t)).\endalign$$
That is, $\f{f'(t)}{f(t)}=\f{x-(2r+1)t}{1+t^2}.$ Solving this
differential equation yields (3.11). $\q\square$
 \pro{Corollary 3.2} For $n\in\Bbb N_0$,
$$D_n^{(r)}(-x)=(-1)^nD_n^{(r)}(x)\qtq{and}D_n^{(r)}(0)=\cases 0&\t{if $n$ is odd,}
\\n!\b{-r-1/2}{n/2}&\t{if $n$ is even.}
\endcases\tag 3.12$$
\endpro
Proof. By Theorem 3.4,
$$ \sum_{n=0}^{\infty}D_n^{(r)}(-x)\f{(-t)^n}{n!}
=(1+t^2)^{-r-\f 12} \;\t{\rm e}^{-x\;arctan\; (-t)} =(1+t^2)^{-r-\f
12} \;\t{\rm e}^{x\;arctan\;
t}=\sum_{n=0}^{\infty}D_n^{(r)}(x)\f{t^n}{n!}.$$ Thus,
$(-1)^nD_n^{(r)}(-x)=D_n^{(r)}(x)$. Taking $x=0$ in Theorem 3.4 and
then applying Newton's binomial theorem we see that
$\sum_{n=0}^{\infty}D_n^{(r)}(0)\f{t^n}{n!}=(1+t^2)^{-r-\f
12}=\sum_{k=0}^{\infty}\b{-r-\f 12}kt^{2k}.$ Comparing the
coefficients of $t^n$ on both sides yields the remaining
part.
$\q\square$
  \pro{Theorem 3.5} For $n\in\Bbb N_0$
we have
$$\align &D_n^{(r)}(x)
=x^n-\sum_{k=1}^{n-1}k(k+2r)D_{k-1}^{(r)}(x) x^{n-1-k}, \tag 3.13
\\&n!d_n^{(r)}(x)=(1+2x)^n+\sum_{k=1}^{n-1}(k+2r)\cdot
k!d_{k-1}^{(r)}(x)(1+2x)^{n-1-k}.\tag 3.14
\endalign$$
\endpro
Proof. For $x\not=0$ and $k=0,1,2,\ldots$ we have
$\f{D_{k+1}^{(r)}(x)}{x^{k+1}}-\f{D_k^{(r)}(x)}{x^k}=-k(k+2r)
\f{D_{k-1}^{(r)}(x)}{x^{k+1}}.$ Thus,
$$\align -\sum_{k=1}^{n-1}k(k+2r)\f{D_{k-1}^{(r)}(x)}{x^{k+1}}
&=\sum_{k=1}^{n-1}\Big(\f{D_{k+1}^{(r)}(x)}{x^{k+1}}
-\f{D_k^{(r)}(x)}{x^k}\Big) =
\f{D_n^{(r)}(x)}{x^n}-\f{D_1^{(r)}(x)}x.\endalign$$ Multiplying by
$x^{n}$ on both sides and noting that $D_1^{(r)}(x)=x$ we deduce
(3.13) for $x\not=0$. When $x=0$, (3.13) is also true by (3.2).
\par By Lemma 3.1, $(-i)^nn!d_n^{(r)}(x)=D_n^{(r)}(-i(1+2x))$. Thus,
$$\align &(-i)^nn!d_n^{(r)}(x)
\\&=D_n^{(r)}(-i(1+2x))=(-i(1+2x))^n-\sum_{k=1}^{n-1}k(k+2r)
D_{k-1}^{(r)}(-i(1+2x))(-i(1+2x))^{n-1-k}
\\&=(-i(1+2x))^n-\sum_{k=1}^{n-1}k(k+2r)
(-i)^{k-1}(k-1)!d_{k-1}^{(r)}(x)(-i(1+2x))^{n-1-k}
\\&=(-i)^n\Big\{(1+2x)^n+\sum_{k=1}^{n-1}(k+2r)\cdot
k!d_{k-1}^{(r)}(x)(1+2x)^{n-1-k}\Big\}.\endalign$$
 This proves (3.14). $\q\square$

 \pro{Corollary 3.3} Let $n\in\Bbb N$. Then
$$\align&[x^n]d_n^{(r)}(x)=\f{2^n}{n!},
\ [x^{n-1}]d_n^{(r)}(x)=\f{2^{n-1}} {(n-1)!},\q
[x^{n-2}]d_n^{(r)}(x)=\f{2^{n-2}}{(n-2)!} \Big(r+\f{n+1}3\Big)\
(n\ge 2),\\& [x^n]D_n^{(r)}(x)=1\qtq{and} [x^{n-2}]D_n^{(r)}(x)=
-\f{(n-1)n(2n-1+6r)}6 \ (n\ge 2),\endalign$$ where $[x^k]f(x)$ is
the coefficient of $x^k$ in the power series expansion of $f(x)$.
\endpro
Proof. From Theorem 3.5 we see that
 $[x^n]D_n^{(r)}(x)=1$ and so
$$[x^{n-2}]D_n^{(r)}(x)=-\sum_{k=1}^{n-1}k(k+2r)=-\sum_{k=1}^{n-1}k^2
-2r \sum_{k=1}^{n-1}k=-\f{(n-1)n(2n-1)}6-rn(n-1).$$ By Theorem 3.5,
$[x^n]d_n^{(r)}(x)=[x^n]\f{(1+2x)^n}{n!}=\f{2^n}{n!}$,
$[x^{n-1}]d_n^{(r)}(x)=[x^{n-1}]\f{(1+2x)^n}{n!}=\f{2^{n-1}}
{(n-1)!}$ and
$$[x^{n-2}]n!d_n^{(r)}(x)=\b n2 2^{n-2}+\sum_{k=1}^{n-1}
(k+2r)k\cdot 2^{k-1}\cdot 2^{n-1-k}=2^{n-2}n(n-1)
\Big(r+\f{n+1}3\Big)\ (n\ge 2).$$ This yields the result.
$\q\square$

 \pro{Theorem 3.6} For any nonnegative integer $n$
we have
$$D_n^{(r)}(x)=D_n^{(r+1)}(x)+n(n-1)D_{n-2}^{(r+1)}(x)=\sum_{k=0}^{[n/2]}\b
n{2k}\b{-r}k(2k)!D_{n-2k}^{(0)}(x).\tag 3.15$$
\endpro
Proof. By Theorem 3.4, for $|t|<1$,
$$\sum_{n=0}^{\infty}D_n^{(r)}(x)\f{t^n}{n!}=(1+t^2)\sum_{n=0}^{\infty}D_n^{(r+1)}(x)\f{t^n}{n!}
=(1+t^2)^{-r}\sum_{n=0}^{\infty}D_n^{(0)}(x)\f{t^n}{n!}.$$ Now
comparing the coefficients of $t^n$ on both sides
yields the result.
$\q\square$
\par
 Finally we state the linearization formula for
$D_m^{(r)}(x)D_n^{(r)}(x)$.
 \pro{Theorem 3.7} Let $m$ and $n$ be nonnegative
integers. Then
$$D_m^{(r)}(x)D_n^{(r)}(x)=\sum_{k=0}^{\min\{m,n\}}\b mk
\b nkk!^2\b{2r+m+n-k}k D_{m+n-2k}^{(r)}(x).\tag 3.16$$
\endpro
Proof. This is immediate from Theorem 2.7 and Lemma 3.1.
$\q\square$
\par\q
\par {\bf Acknowledgments}
\par The author was supported by the
National Natural Science Foundation of China (Grant No. 11371163).

\end{document}